
\input amstex
\documentstyle{amsppt}

\magnification=\magstep 1

\def \rank {\operatorname{rank}}
\def \codim {\operatorname{codim}}
\def \point {\operatorname{point}}
\def \Spec {\operatorname{Spec}}
\def \dim {\operatorname{dim}}
\def \top {\operatorname{top}}
\def \Ker {\operatorname{Ker}}
\def \C {\operatorname{C}}
\def \det{\operatorname{det}}
\def \card {\operatorname{card}}
\def \deg {\operatorname{deg}}
\def \graph {\operatorname{graph}}
\def \Pf {\operatorname{Pf}}
\def \Coker {\operatorname{Coker}}

\let\toto=\twoheadrightarrow
\let\hoto=\hookrightarrow
\let\loto=\longrightarrow

\let\ge=\geqslant
\let\le=\leqslant

\document
\topmatter
\title
SCHUR $Q$-FUNCTIONS AND
DEGENERACY LOCUS FORMULAS FOR MORPHISMS WITH SYMMETRIES
\endtitle
\rightheadtext{LASCOUX AND PRAGACZ}
\leftheadtext{DEGENERACY LOCUS FORMULAS}
\author
Alain Lascoux  and  Piotr Pragacz
\endauthor
\subjclass 14M12, 14C17, 14M15
\endsubjclass
\thanks 
This research was supported by the grant No.5031 of French-Polish
cooperation C.N.R.S. -- P.A.N., and the KBN grant No.2P03A 05112.
\endthanks

\endtopmatter

\vskip30pt
\rightline{{\it To Maria}}
\vskip15pt

\noindent
{\bf Abstract.}
{\eightrm We give closed-form formulas for the fundamental classes
of degeneracy loci associated with vector bundle maps given locally
by (not necessary square) matrices which are symmetric 
(resp. skew-symmetric) w.r.t. the
main diagonal. Our description uses essentially Schur $Q$-polynomials
of a bundle and is based on a push-forward formula for these
polynomials in a Grassmann bundle, established in [P4].}

\vskip40pt

\rightline{\eightrm ``Something which is not testable is not 
scientific."}
\rightline{\eightrm  Poper's criterium}

\vskip5pt

\head
{\bf Introduction}
\endhead
 
The goal of the present paper is to state and prove new  
closed-form formulas for the fundamental classes of some degeneracy loci.

We will be here interested
in degeneracy loci associated with vector bundle maps given locally
by (not necessary square) matrices that are symmetric 
(resp. skew-symmetric) w.r.t. the main diagonal. To be more precise, 
let $\alpha: E\toto F$ be
a surjection of two vector bundles of respective ranks $e$ and $f$
on a variety $X$. 
We denote by $E\vee F$ (resp. $E\wedge F$)
the kernel of the surjection 
$$
\CD
E\otimes F @> \alpha \otimes 1 >> F\otimes F \toto \wedge^2F
\endCD
$$ 
(resp. 
$$
\CD
E\otimes F @> \alpha \otimes 1 >> F\otimes F \toto S^2F).
\endCD
$$
In the present paper, we will call a morphism
$\varphi: E^* \to F$ {\it symmetric} (resp. {\it skew-symmetric}) provided
it is induced by a section of the subbundle
$E\vee F$ (resp. $E\wedge F$) of $E\otimes F$. (Note that for $E=F$
these notions coincide with the usual notions of ``symmetric" 
and ``skew-symmetric" morphisms.)

\smallskip

Our goal is to describe the fundamental classes of the loci 
$$
D_r(\varphi)=\{x\in X: \rank \ \varphi(x)\le r\}
$$ 
with the help of some explicitly given polynomials in the Chern classes of 
$E$ and $F$. Set $n=e-f$. More precisely, using the terminology of [F-P], 
the degeneracy
locus $D_r(\varphi)$ is represented, in the symmetric case, by the polynomial
$$
\sum_I Q_{(f-r,f-r-1,...,1)+I}(F)\cdot s_{\C\widetilde I}(E-F),
\tag *
$$
where the sum is over all partitions $I$ in $(n)^{f-r}$, and for the conjugate
partition $\widetilde I=(\widetilde i_1,\ldots,\widetilde i_n)
\subset (f-r)^n$ we write $\C\widetilde I 
= (f-r-\widetilde i_n,\ldots, f-r-\widetilde i_1)$ for the ``complement"
of $\widetilde I$ in $(f-r)^n$.
Here, $Q_J({\sssize\bullet})$ 
are Schur $Q$-polynomials whose definition is recalled in Section 1.
The connection of Schur $Q$-polynomials to geometry was originally
established in [P2] in two contexts. Firstly, in [P2, Sect.7], the ideal
of polynomials universally supported on the $r$th symmetric degeneracy locus
was described in terms of Schur $Q$-polynomials. This is closely connected
with the present paper, see below. Secondly, in [P2, (8.7)] 
the connection of Schur $Q$-functions to the cohomology rings 
of isotropic Grassmannians was given. 
(This subject was then developed in [P3, Sect.6].)
Consult also [P4] and [F-P] for more about geometric properties 
of $Q$-polynomials, and especially
for their applications to intersection theory and enumerative geometry.

By virtue of [P2, Sect.7], the form of polynomials (*)
representing $[D_r(\varphi)]$ is not at all surprising.
Observe that
for $\varphi'=\varphi \circ \alpha^* : F^* \to E^* \to F$, we have
$D_r(\varphi) \subset D_r(\varphi')$, so the polynomials representing
$D_r(\varphi)$ are universally supported on the $r$th (symmetric)
degeneracy locus associated with symmetric bundle maps from rank $f$
vector bundles to theirs duals, in the sense of [P2, Sect.7] (see 
also [F-P, Sect. 4.4]). Therefore, by [P2, Th.7.2 and Prop.7.17], 
the polynomial representing $[D_r(\varphi)]$ is a $\Bbb Z[c.(E),c.(F)]$-
combination of $Q_{(f-r,f-r-1,\ldots,1)+I}(F)$ for $I$ contained
in $(r)^{f-r}$. For $n\le r$, this is exactly the combination (*).
For $n > r$, this suggests that at the cost of replacing in (*) 
the coefficients from $\Bbb Z[c.(E-F)]$ 
by the coefficients from a larger ring
$\Bbb Z[c.(E),c.(F)]$, one can make the sum (*) smaller. We also note
that to get (*), we use exactly the same instance of the Gysin map formula,
recalled in (1.8) below, as to get the just mentioned results from
[P2, Sect.7]. 
\smallskip

Similar formulas and discussion hold true for the degeneracy loci
associated with skew-symmetric morphisms, see Section 3.

\smallskip

When $E=F$, our formulas specialize to the ones given by 
J\'ozefiak and the authors in [J-L-P], and Harris and Tu
in [H-T] (see also [P1]). 
Recall that all the
three last mentioned papers gave a ``modern treatment" \`a la ``Thom-Porteous"
of the formulas for the degree of determinantal varieties in the spaces of
symmetric and skew-symmetric matrices of forms, 
given classically by Giambelli in [G]. Some special cases of the formulas 
from [H-T], [J-L-P] and [P1] were given independently by Barth [Ba], Damon 
and Tyurin.

\smallskip

When $E=\bigoplus^e_{i=1}\Cal O(p_i)\,, \ F=\bigoplus^f_{i=1}\Cal O(p_i)$ 
are vector bundles over a projective space, our formulas give the degrees
of determinantal varieties defined by matrices of forms, satysfying
the above symmetry conditions. This problem was considered, in the symmetric
case, 
by Bottasso in [Bo], generalizing Giambelli's study to the case of
not necessary square matrices. 
Bottaso did not give closed-form formulas
for these degrees, but established some recursions for their computations.
The present paper offers a modern treatment and closed-form version 
of [Bo].

\smallskip

There are several motivations to do such computations in intersection
theory. These computations combine intersection theory on homogeneous spaces
(notably Grassmannians and flag varieties) and algebra of different
types of symmetric functions. The obtained formulas have applications
to the enumerative theory of singularities, Brill-Noether theory etc.
Again, we refer the reader to [P4], and to [F-P] for more information 
on these matters.

\smallskip

We give in the present paper formulas for the ``Thom-Porteous" case.
It would be interesting to generalize them to the ``quiver variety"
case, similarly to the work of Buch and Fulton [B-F].

\smallskip

Since the degeneracy loci studied here are intimately connected with 
some Schubert varieties in symplectic and orthogonal Grassmannians
(see Section 1),
the formulas given in the present paper should be gotten ``in principle" 
from the divided
difference operator approach developed in [F2-3] from one side, and in 
[P-R] and [L-P] from another one. 
We do not see, however, how to do it conceptually using divided
differences. 
\medskip

The article is organized as follows. 

\smallskip

In Section 1, we recall the definitions
and properties of two families of symmetric polynomials that we need
in the present paper: Schur polynomials in a difference of bundles
and Schur $Q$-polynomials. We also recall
a formula for Gysin push-forward of Schur $Q$-polynomials
in a Grassmann bundle from [P4], which is basic to the present paper. 
Usually, we do not state results in their full generality, but only
in the range needed in this paper.
Moreover, we give some preliminary properties and examples 
of the degeneracy loci associated with symmetric and skew-symmetric 
morphisms.

In Section 2, we give some formulas for the top Chern classes of the
bundles $E\vee F$ and $E\wedge F$, generalizing the formulas from [L].

In Section 3, we prove the main formula of the present paper. The proof
follows a well-known pattern: it
uses essentially a certain desingularization of $D_r(\varphi)$ for a 
``universal" $\varphi$ and the above mentioned push-forward formula for 
Schur $Q$-polynomials. Some examples are also discussed.

In Section 4, we discuss some variations of computation of the polynomials
representing $D_r(\varphi)$, based on the technique of ``constructions 
with a nontrivial generic fibre'' invented in [P1, Sect.2]. 
This material is recalled in Proposition 4.1, where, in fact, a certain 
straightening  of this last-mentioned method is presented. 
As an application, we get some algebraic equations involving
Gysin maps on their LHS's and some closed-form expressions involving
symmetric polynomials on their RHS's.

\bigskip

\head
{\bf 1. Preliminaries: recollection of some definitions and results}
\endhead
\smallskip

In the present paper, we will need two families of symmetric polynomials.
We recall now their definitions and needed properties.
Let $A=(a_1,\ldots,a_n)$ and $B=(b_1,\ldots,b_m)$ be two
sequences of commuting elements in a ring.

\smallskip
\noindent
{\bf (1.1)} \  Let $I=(i_1,\dots,i_k)$  be a partition.
Recall that by a  {\it partition}  (of some natural number) we understand 
a sequence of integers  
$I=(i_1,\ldots,i_k)$, where  $i_1\ge i_2\ge \ldots\ge i_k\ge 0$ and the sum 
$\sum i_p$ is the partitioned number.  
We define 
$$
s_I(A-B) = \det \Bigl[s_{i_p-p+q}(A-B)\Bigr]_{_{1\le p,q\le k}} ,
$$
where  $s_i(A-B)$  is defined by
$$
\sum\limits^{\infty}_{i=-\infty} s_i(A-B)=
\prod\limits^n_{p=1}(1-a_p)^{-1}\ \prod\limits^m_{q=1}(1-b_q) .
$$
Moreover, we put $s_I(A)=s_I(A-B)$ for $B=(0,\ldots,0)$. When $A$ and $B$
are sequences of variables, then
the polynomials $s_I(A)$ are called {\it Schur polynomials}
or  $S$-{\it polynomials} and
the polynomials $s_I(A-B)$ are often called  
{\it Schur polynomials in a difference of alphabets}. 
For example, the ``resultant" can be expressed in a closed form as
$$
\prod_{i,j}(a_i-b_j)=s_{(m,\ldots,m)}(A-B),
\tag 1.1.1
$$
where the part $m$ appears $n$ times.
\smallskip
\noindent
{\bf (1.2)} \ Let $Q_i(A)$ be defined by the expansion
$$
\sum\limits^{\infty}_{i=-\infty} Q_i(A) =
\prod\limits^n_{p=1}(1+a_p)(1-a_p)^{-1}. 
$$  
Given positive integers $i>j$, we set
$$
Q_{(i,j)}(A)=Q_i(A)\cdot Q_j(A) + 2\sum\limits^j_{p=1}(-1)^p Q_{i+p}(A)\cdot
 Q_{j-p}(A).
$$
Finally, if $I=(i_1 > \ldots > i_k > 0)$  is a {\it strict} partition,
then for odd $k$ we put
$$
Q_I(A) = \sum_{p=1}^k (-1)^{p-1}Q_{i_p}(A)\cdot Q_{(i_1,\ldots,i_{p-1},i_{p+1},
\ldots,i_k)}(A),
$$
and for even $k$,
$$
Q_I(A) = \sum_{p=2}^k (-1)^pQ_{(i_1,i_p)}(A)\cdot Q_{(i_2,\ldots,i_{p-1},
i_{p+1},\ldots,i_k)}(A).
$$
These define the $Q_I(A)$'s by recurrence on $k$.
When $A$ is a sequence of variables, then the polynomials $Q_I(A)$ 
are called  {\it Schur Q-polynomials} after [Sch].
For example, one has
$$
Q_{(n,n-1,\ldots,1)}(A)=\prod_{i\le j}(a_i+a_j)=2^n s_{(n,n-1,\ldots,1)}(A).
\tag 1.2.1
$$

\smallskip
\noindent
{\bf (1.3)} \ Let $E$ and $F$ be two vector bundles. 
Then  $s_I(E-F)$ is  defined to be $s_I(A-B)$, where $A$ and $B$ are
the sequences of the Chern roots of $E$ and $F$ respectively.  
Similarly, set $s_I(E)=s_I(A)$ \ and $Q_I(E)=Q_I(A)$.

\smallskip
\noindent
{\bf (1.4)}
For partitions $I$ and $J$,  we write  
$I\supset J$  if $i_1\ge j_1$, $i_2\ge j_2$, ... . Moreover,
we denote by 
$I+J$ the partition $(i_1+j_1,i_2+j_2,\ldots)$.
 
The ``rectangular" partition $(i,\dots,i)$ ($r$-times) is denoted by $(i)^r$ 
and the ``triangular" partition $(k,k-1,\dots,2,1)$ is denoted by $\rho_k$. 

For a given partition $I$, $l(I)= \card \{ p :i_p > 0 \}$ denotes its length 
and $\widetilde I$ denotes the partition conjugate of $I$, i.e., 
$\widetilde I=(\widetilde i_1,\widetilde i_2,\dots)$ where $\widetilde i_p
=\card \{ q:i_q\ge p \}$.

As is common, we will often omit brackets writing partitions in lower
indices.

\smallskip
\noindent
{\bf (1.5)} \ Recall  that for every strict partition 
$I=(i_1>\dots>i_k>0)$ and $A$ a sequence of variables, one has 
$$
Q_I(A)= 2^k P_I(A)
$$
for some polynomial $P_I(A)$ with integer coefficients. 
These polynomials are called {\it Schur P-polynomials}. 
For example,
$$
P_{\rho_{n-1}}(A)=\prod_{i<j}(a_i+a_j)=s_{\rho_{n-1}}(A).
$$

Given a vector
bundle $E$, we set $P_I(E)=P_I(A)$, where $A$ is specialized to the sequence
of the Chern roots of $E$. 
For another expression of $P_I(A)$ in the form of a quadratic polynomial 
in the $s_J(A)$'s, see [L-L-T].
\smallskip
\noindent
{\bf (1.6)} \ We recall the following two {\it factorization formulas}. 
Let $I$ be a partition such that $l(I)\le n$. Then  
$$
s_{(m)^n+I}(A-B)= s_{(m)^n}(A-B)\cdot s_I(A)
\tag 1.6.1
$$
and
$$
Q_{\rho_{n-1}+I}(A)=Q_{\rho_{n-1}}(A)\cdot s_I(A).
\tag 1.6.2
$$
We refer to [B-R] [L-S] and [St] (see also [M], [P3] and [L-L-T])
for more about these formulas. In particular, the proof of (1.6.1)
given in [L-S] is valid for any elements $A$ and $B$ of ranks $n$
and $m$ in a $\lambda$-ring.
\medskip
\noindent
{\bf 1.7} \ We recall the following formulas for the top Chern classes
of some tensor operations. We write, in the present paper,
$c_{\top}(E)$ for the top Chern
class of a bundle $E$. Let $E$ and $F$ be two vector bundles of ranks
$e$ and $f$ respectively. We have
$$
c_{\top}(E\otimes F)=\sum_I s_I(E)\cdot s_{\C\widetilde I}(F),
\tag 1.7.1
$$
where the sum is over all partitions in $(f)^e$ and for the conjugate
partition $\widetilde I \subset (e)^f$ \ we write $\C \widetilde I=
(e-\widetilde i_f,\ldots, e-\widetilde i_1)$ for the ``complement"
of $\widetilde I$ in $(e)^f$. Also,
$$
c_{\top}(S^2E) = 2^e s_{\rho_e}(E) = Q_{\rho_e}(E),
\tag 1.7.2
$$
and
$$
c_{\top}(\wedge^2E) = s_{\rho_{e-1}}(E) = P_{\rho_{e-1}}(E).
\tag 1.7.3
$$
We refer to [L] (see also [M, p.47-48 and p.67]) for more details.
\smallskip
\noindent
{\bf 1.8} \ We recall the following push-forward formula from 
[P4]. Let  $\pi: G=G^q(E) \to X$ be 
the Grassmann bundle parametrizing  $q$-quotients of a vector bundle $E$
on a variety $X$. 
Write $r=e-q$. Let
$$
0 \to R \to E_G \to Q \to 0
$$
be the tautological sequence on $G$ with $\rank R = r$ and $\rank Q = q$. 
Let $I=(i_1>\ldots>i_k>0)$ be a  
strict partition with $k\le q$. Then for $\alpha \in
A_*(X)$, we have 
$$
\pi_*\Bigl[c_{\top}(R\otimes Q)\cdot P_I(Q) \cap \pi^*\alpha \Bigr] 
= d\cdot P_I(E)\cap \alpha, 
$$
where $d$ is zero if $(q-k)r$ is odd, \ and 
$$
d = {[(e-k)/2]\choose[(q-k)/2]}
$$ 
in the opposite case. Here, the symbol $[ \ {\sssize\bullet} \ ]$ 
means the integer
part of a rational number. For a proof, we refer to [P4, App.1]. 
\smallskip
 
We will need three special instances of this formula. 
First, suppose that $I$ is a strict partition with $l(I)=q$.
Then 
$$
\pi_*\Bigl[c_{\top}(R \otimes Q)\cdot Q_I(Q) \cap \pi^*\alpha \Bigr]= 
Q_I(E) \cap \alpha. 
\tag 1.8.1
$$
Secondly, assume that $I$ is a strict partition with $l(I)=q-1$ 
and $r$ is even. Then
$$
\pi_*\Bigl[c_{\top}(R \otimes Q)\cdot P_I(Q) \cap \pi^*\alpha \Bigr]= 
P_I(E) \cap \alpha. 
\tag 1.8.2
$$
Thirdly, suppose that $I$ is a strict partition with $l(I)=q-1$ and
$r$ is odd. Then
$$
\pi_*\Bigl[c_{\top}(R \otimes Q)\cdot P_I(Q) \cap \pi^*\alpha \Bigr]=0.
\tag 1.8.3
$$
\medskip
\noindent
{\bf 1.9} \ Recall that the Gysin map in a Grassmann bundle admits the
following explicit description. Let $(a_1,\ldots,a_q)$ be the sequence of the
Chern roots of $Q$ and $(a_{q+1},\ldots,a_e)$ be the sequence of the Chern
roots of $R$. Then, writing $A=(a_1,\ldots,a_e)$ for the sequence of the
Chern roots of $E$, the Gysin map in question is induced by the following 
symmetrizing operator. 
Let $S_e$ be the group of permutations of $(1,\ldots,e)$, $S_q$ the
group of permutations of $(1,\ldots,q)$, and $S_r$ the group of
permutations of $(q+1,\ldots,e)$.
For $P\in\Bbb Z[A]^{S_q\times S_r}$, the symmetrizing operator in question
acts as follows:
$$
P \longmapsto \sum\limits_{\mathstrut\overline{\sigma}\in S_e/S_q\times S_r}
\sigma\left({P \big/ \prod\limits_{i \le q<j}(a_i-a_j)}
\right).
$$
For more on this, see, e.g., [P4, Sect.4].
\smallskip
\noindent
{\bf 1.10} \ We now switch to the setup of the Introduction. 
In the present paper, to be on the safe side, we assume
that the ground field $k$ is algebraically closed of characteristic
different from $2$.  
This is because in Section 1 and 4 we make
use of isotropic symplectic and orthogonal Grassmannians. 
Let $U\toto V$ be vector spaces of dimensions $e$ and $f$, respectively. 
Let $X$ denote the affine space $U\vee V$ (resp. $U\wedge V$). 
In this situation, there exists a tautological morphism 
$\varphi: (E=U_X)^* \to (F=V_X)$. 
For this $\varphi$, $D_r(\varphi)$ is the restriction to an appropriate
open set of some Schubert variety in the 
symplectic (resp. orthogonal) Grassmannian of $f$-dimensional isotropic
subspaces in $k^{2e}$. 
More precisely, let $\beta: U^* \to V$ be a linear map and $\gamma: V^* \to U$
its dual. Consider $U^*\oplus U$ equipped with canonical nondegenerate 
bilinear forms
$$
\langle(v_1,u_1),(v_2,u_2)\rangle = v_1(u_2) \pm v_2(u_1)  
\ \ \ \ u_i\in U, v_i\in U^*,
$$
with the + sign giving a symmetric form and -- sign a skew-symmetric form.
The assignment \ $\beta \mapsto \graph (\gamma)$ \ embeds $U\vee V$
(resp. $U\wedge V$) as an open subset $A$ of the Grassmannian $G$ of
$f$-dimensional isotropic subspaces of $U^*\oplus U$ w.r.t. the just
defined skew-symmetric (resp. symmetric) form on $U^*\oplus U$. Then
$D_r(\varphi)$ is the restriction to $A$ of the ``determinantal" Schubert 
subvariety of $G$ parametrizing those $f$-dimensional isotropic
subspaces of $U^*\oplus U$, which intersect the maximal isotropic
subspace $U^*\oplus 0$ in dimension $\ge f-r$. (When $\varphi$ is 
skew-symmetric and $e=f$, we assume $r$ to be even.)
Hence $D_r(\varphi)$ is irreducible, normal and Cohen-Macaulay (by results
of De Concini and Lakshmibai [DC-L]); moreover its
codimension $c(r)$ equals 
$$
(e-f)(f-r)+(f-r)(f-r+1)/2 \ \ \  (\hbox{resp.} \ (e-f)(f-r)+(f-r)(f-r-1)/2).
$$
Perhaps the easiest way to remember this number is to set $n=e-f$ and $q=f-r$
(we will keep this notation throughout the rest of the paper),
and note that $c(r)$ equals
$$
nq+q(q+1)/2=q(2n+q+1)/2 \ \ \  (\hbox{resp.} \ nq+q(q-1)/2=q(2n+q-1)/2).
$$
In general, for $\varphi: E^* \to F$ as in the Introduction, 
the scheme structure on $D_r(\varphi)$ is defined as follows.
Set $\tilde{X}=\Spec \ S^{\sssize\bullet}(E\vee F)^*$ (resp.
$\tilde{X}= \Spec \ S^{\sssize\bullet}(E\wedge F)^*$\,).
Observe that $\varphi$ induces a section $s:X\to \tilde{X}$\,.
There exists the tautological bundle homomorphism $\tilde{\varphi}:
\tilde E^*\to\tilde F$ where $\tilde E=E_{\tilde{X}}$\,,
$\tilde F=F_{\tilde{X}}$\, \ such that $s^*(\tilde{\varphi})=\varphi$\,.
Then, defining first the scheme structure on $D_r(\tilde{\varphi})$ as
above with the help of an appropriate Schubert bundle in an isotropic
Grassmann bundle, we define the scheme structure on $D_r(\varphi)$
as the schematic preimage via $s$ of the one on $D_r(\tilde{\varphi})$.

To the best of our knowledge, the above determinantal varieties
have not been studied algebraically 
like the ``ordinary" determinantal varieties were studied by
Eagon and Hochster or using the Hodge algebra technique.
Probably such a study would allow us to formulate the main results
of the present paper with less restrictive assumptions on the ground field.
\medskip
\noindent
{\bf 1.11} \ We pass now to some examples. 
\smallskip
\noindent
(1.11.1) \ Let $\varphi: E^* \to F$ be symmetric, where $e\ge f$ are
arbitrary and $r=f-1$. Then the expected codimension of $D_r(\varphi)$
is $e-f+1$, so we are in the situation of the Giambelli-Thom-Porteous
formula for the locus defined by maximal minors.
We infer that $D_r(\varphi)$ is represented by $s_{e-f+1}(F-E^*)$.
For example, for $e=4, f=3$, writing formally by the splitting principle
$c(E)=c(F)(1+d)$, we get
$$\aligned
&s_2(F-E^*)=s_2(F)-s_1(F)s_1(E^*)+s_{1,1}(E^*)=s_2(F)+s_1(F)s_1(E)+s_{1,1}(E)
\cr
&=s_2(F)+s_1(F)\bigl(s_1(F)+d\bigr)+s_{1,1}(F)+s_1(F) d
=2\bigl(s_2(F)+s_{1,1}(F)+s_1(F) d\bigr) \cr
&=2\bigl(s_2(F)+s_{1,1}(E)\bigr)=2s_1(E)s_1(F).
\endaligned
$$
For $e=5, f=3$, we get the representing polynomial
$$
s_3(F-E^*)=2\bigl(s_1(E)s_2(F)+s_{1,1,1}(E)\bigr)
=2\bigl(s_3(F)+s_{1,1}(E)s_1(F)\bigr).
$$
\medskip
\noindent
(1.11.2) \ In this example, $X$ is a projective space, $F=\Cal O(a) \oplus
\Cal O(b) \oplus \Cal O(c)$ and $r=2$.

Assume first that $E=\Cal O(a) \oplus \Cal O(b) \oplus \Cal O(c) \oplus 
\Cal O(d)$ and $\varphi: E^* \to F$ is given by a matrix of forms

$$\left(\ \CD
0 &\quad a_{12} &\quad a_{13} &\quad x \\ 
-a_{12} &\quad 0 &\quad a_{23} &\quad y \\
-a_{13} &\quad -a_{23} & \quad 0 &\quad z
\endCD\ \right).$$
\smallskip
\noindent
Then $D_2(\varphi)$ is of codimension $1$, defined by
$$
a_{12}z - a_{13}y + a_{23}x = 0.
$$
The degree of this equation is $a+b+c+d=s_1(E)$.

\medskip

Assume now that $E=\Cal O(a) \oplus \Cal O(b) \oplus \Cal O(c) \oplus 
\Cal O(d) \oplus \Cal O(e)$ and $\varphi : E^* \to F$ is given by a matrix
of forms

$$\left(\ \CD
0 &\quad a_{12} &\quad a_{13} &\quad x &\quad w \\ 
-a_{12} &\quad 0 &\quad a_{23} &\quad y &\quad v \\
-a_{13} &\quad -a_{23} & \quad 0 &\quad z &\quad u 
\endCD\ \right).$$
\smallskip
\noindent
Then $D_2(\varphi)$ is of codimension $2$, defined by
$$
a_{12}z - a_{13}y + a_{23}x = 0  \ \ \ \hbox{and} \ \ \ 
a_{12}u - a_{13}v + a_{23}w = 0.
$$
By B\'ezout's theorem, the degree of $D_2(\varphi)$ is
$(a+b+c+d)(a+b+c+e)=s_2(F) + s_{1,1}(E)$.

\medskip
\noindent
(1.11.3) \ Let $\varphi: E^*\to F$ be skew-symmetric, $n=e-f=1$ and
$r$ is an arbitrary even nonnegative number less than $f$. 
Assume that $\varphi$ is given locally by a matrix
\smallskip
$$\left(\ \CD
0 &\quad a_{12} &\quad .\;.\;. \quad & a_{1f} &\quad  b_1 \\
-a_{12} &\quad 0 &\quad .\;.\;. \quad & a_{2f} &\quad  b_2 \\
\vbox {\offinterlineskip\hbox{.}\vskip 3pt\hbox{.}\vskip 3pt\hbox{.}}&\quad
\vbox {\offinterlineskip\hbox{.}\vskip 3pt\hbox{.}\vskip 3pt\hbox{.}}&
\vbox{\offinterlineskip\hbox{.}\vskip 3pt \hbox {\;\;.}\vskip 3pt
\hbox {\;\;\;\;.}}&
\vbox {\offinterlineskip\hbox{.}\vskip 3pt\hbox{.}\vskip 3pt\hbox{.}}&\quad
\vbox {\offinterlineskip\hbox{.}\vskip 3pt\hbox{.}\vskip 3pt\hbox{.}} \\
-a_{1f} &\quad -a_{2f}&\quad .\;.\;. \quad & 0 &\quad  b_f
\endCD \ \right).$$
Then $D_r(\varphi)$ is locally defined by the $(r+2)$-Pfaffians of the extended
$e\times e$ skew-symmetric matrix
\smallskip
$$\left(\ \CD
0 &\quad a_{12} &\quad .\;.\;. \quad & a_{1f} &\quad  b_1 \\
-a_{12} &\quad 0 &\quad .\;.\;. \quad & a_{2f} & \quad  b_2 \\
\vbox {\offinterlineskip\hbox{.}\vskip 3pt\hbox{.}\vskip 3pt\hbox{.}}&\quad
\vbox {\offinterlineskip\hbox{.}\vskip 3pt\hbox{.}\vskip 3pt\hbox{.}}&
\vbox{\offinterlineskip\hbox{.}\vskip 3pt \hbox {\;\;.}\vskip 3pt
\hbox {\;\;\;\;.}}&
\vbox {\offinterlineskip\hbox{.}\vskip 3pt\hbox{.}\vskip 3pt\hbox{.}}&\quad
\vbox {\offinterlineskip\hbox{.}\vskip 3pt\hbox{.}\vskip 3pt\hbox{.}}\\
-a_{1f} &\quad -a_{2f}&\quad .\;.\;. \quad & 0 &\quad  b_f \\
-b_1 &\quad -b_2 &\quad .\;.\;. \quad & -b_f &\quad  0
\endCD \ \right).$$
This is seen by the following ``universal local study". Assume that the
$a_{ij}$'s and the $b_p$'s are variables over $k$. Let $\Cal J$ be the ideal
generated by $(r+1)$-minors of the former matrix, and let $\Cal I$
(resp. $\Cal P$) be the ideal generated by $(r+1)$-minors (resp. 
$(r+2)$-Pfaffians) of the latter matrix. Of course, $\Cal J\subset \Cal I$;
moreover, $\Cal I \subset \Cal P$ by [B-E, p.462]. Since $\Cal P$ is
a prime ideal of height equal to the expected codimension $c(r)$
(see [K-L]), the ideal $\Cal P$ defines (locally) the scheme structure
on $D_r(\varphi)$. Hence using [G], ..., we see that $D_r(\varphi)$
is represented by \ $s_{\rho_{e-r-1}}(E)$. 
For example, for $e=5, f=4, r=2$, the representing polynomial 
equals \ $s_{2,1}(E)$.
\medskip

In (1.11.4-5) we assume that the ambient space $X$ is a projective space.

\smallskip
\noindent
(1.11.4) \ Suppose that $E=\Cal O(a)\oplus \Cal O(b) \oplus \Cal O(c)$, 
$F=\Cal O(a)\oplus \Cal O(b)$, and $\varphi: E^* \to F$ is given by 
a matrix of forms
$$\left(\ \CD
0 &\quad a_{12} &\quad x \\ 
-a_{12} &\quad 0 &\quad y
\endCD\ \right).$$
\smallskip
\noindent
Then $D_1(\varphi)$ is of codimension $1$, defined by the equation
$a_{12}=0$ whose degree is $a+b=s_1(F)$.

More generally, assume that $f$ is even, $E=\Cal O(p_1)\oplus \ldots 
\oplus \Cal O(p_{f+1})$, $F=\Cal O(p_1)\oplus \ldots \oplus \Cal O(p_f)$, 
and $\varphi : E^* \to F$ is given by a matrix of forms
\smallskip
$$\left(\ \CD
0 &\quad a_{12} &\quad .\;.\;. \quad & a_{1f} &\quad  b_1 \\
-a_{12} &\quad 0 &\quad .\;.\;. \quad & a_{2f} &\quad  b_2 \\
\vbox {\offinterlineskip\hbox{.}\vskip 3pt\hbox{.}\vskip 3pt\hbox{.}}&\quad
\vbox {\offinterlineskip\hbox{.}\vskip 3pt\hbox{.}\vskip 3pt\hbox{.}}&
\vbox{\offinterlineskip\hbox{.}\vskip 3pt \hbox {\;\;.}\vskip 3pt
\hbox {\;\;\;\;.}}&
\vbox {\offinterlineskip\hbox{.}\vskip 3pt\hbox{.}\vskip 3pt\hbox{.}}&\quad
\vbox {\offinterlineskip\hbox{.}\vskip 3pt\hbox{.}\vskip 3pt\hbox{.}} \\
-a_{1f} &\quad -a_{2f}&\quad .\;.\;. \quad & 0 &\quad  b_f
\endCD \ \right).$$
Then $D_{f-1}(\varphi)$ is of codimension $1$, defined by the Pfaffian
of the $f\times f$ matrix $A=(a_{ij})$. Indeed, every $(f-1)$-minor
of $A$ occuring in the Laplace expansion along the last column of
an $f$-minor that contains the $(f+1)$th
column, is a multiple of $\Pf(A)$ by already quoted result from [B-E].
The degree of $\Pf(A)$ is $p_1+\ldots+p_f=s_1(F)$.
\medskip
\noindent
(1.11.5) \ Assume that $E=\Cal O(a)\oplus \Cal O(b) \oplus \Cal O(c) \oplus
\Cal O(d)$, $F=\Cal O(a) \oplus O(b)$, and $\varphi: E^* \to F$ is given
by a matrix of forms

$$\left(\ \CD
0 &\quad a_{12} &\quad x &\quad z \\ 
-a_{12} &\quad 0 &\quad y &\quad w 
\endCD\ \right).$$
\smallskip
\noindent
Then $D_1(\varphi)$ is of codimension $2$, defined by
$$
a_{12}=0 \ \ \ \ \ \hbox{and} \ \ \ \ \ xw-zy=0.
$$
By B\'ezout's theorem the degree of $D_1(\varphi)$ is $(a+b)(a+b+c+d)=
s_1(F)s_1(E)$.

\medskip

\head
{\bf 2. Formulas for the top Chern classes of $E\vee F$ and $E\wedge F$}
\endhead
\smallskip

We keep the notation from the Introduction and Section 1. 
In this section, we give
two alternative expressions for the top Chern classes of $E\vee F$ and 
$E\wedge F$. In particular, these expressions give the polynomials
representing $D_r(\varphi)$ for $r=0$. 

\smallskip

In the following, we set $K=\Ker(\alpha: E\toto F)$.

\medskip

\proclaim{\bf Proposition 2.1} \ We have 

\smallskip

\centerline{$c_{\top}(E\vee F)=
\sum_I Q_{\rho_f+I}(F)\cdot s_{\C\widetilde I}(E-F)$,}
\smallskip
\noindent
where the sum is over all partitions in $(n)^{f}$ and for the conjugate
partition $\widetilde I\subset (f)^n$ we write $\C\widetilde I 
= (f-\widetilde i_n,\ldots, f-\widetilde i_1)$. Similarly,

\smallskip

\centerline{$c_{\top}(E\wedge F)=
\sum_I P_{\rho_{f-1}+I}(F)\cdot s_{\C\widetilde I}(E-F)$,}
\smallskip
\noindent
the same sum as above.
\endproclaim
\smallskip
\demo{Proof} \ We give here the proof of the proposition for the bundle 
$E\vee F$, the case $E\wedge F$ being similar. 
We have in the Grothendieck group $K(X)$,
$$
[E\vee F]= [E \otimes F] - [\wedge^2 F] = [S^2F] + [K\otimes F].
$$
Invoking the formulas
$$
c_{\top}(S^2F) = 2^fs_{\rho_f}(F) = Q_{\rho_f}(F),
$$
and
$$
c_{\top}(K\otimes F) = \sum_I s_I(F)\cdot s_{\C\widetilde I}(K),
$$
the same sum as above (see (1.7.1) and (1.7.2)), the assertion follows 
by virtue of the factorization formula for $Q$-polynomials (1.6.2) 
and the equality $[K]=[E]-[F]$ in $K(X)$.
\qed
\enddemo

\medskip

We will see in the next section that the expressions given in the
previous proposition generalize, in a natural way, for higher $r$.
The expressions given in the next proposition do not admit such
a generalization.

\smallskip

\proclaim{\bf Proposition 2.2} \ We have 
\smallskip
\centerline{$c_{\top}(E\vee F) = 2^f \sum_I s_{(e,e-1,\ldots,n+2,n+1)/I}(F)
\cdot s_{\widetilde I}(E-F),$}
\smallskip
\noindent
where the sum runs over all partitions $I\subset (e,e-1,\ldots,n+2,n+1)$. 
Similarly,
\smallskip

\centerline{$c_{\top}(E\wedge F) = \sum_I s_{(e-1,e-2,\ldots,n+1,n)/I}(F) 
\cdot s_{\widetilde I}(E-F),$}
\smallskip
\noindent
where the sum runs over all partitions $I\subset (e-1,e-2,\ldots,n+1,n)$.
\endproclaim
\smallskip
 
\demo{Proof} \ We give here the proof of the proposition for the bundle 
$E\vee F$, the case $E\wedge F$ being similar. 
By the splitting principle, a (more general) question concerning the total 
Chern class of $E\vee F$ leads to the calculation of the product
$$
\prod\limits_{i\le j}(1+a_i+a_j) \prod\limits_{i,j} (1+a_i+b_j),
\tag 2.3
$$
where, formally, 
$$
c(F)=\prod\limits_{i=1}^f(1+a_i)  \ \ \hbox{and} \ \  
c(K)=\prod\limits_{j=1}^n (1+b_j).
$$ 
Write $A=(a_1,\ldots,a_f)$ and 
$B=(b_1,\ldots,b_n)$. Let $a_i^+=1+2a_i$, $1\le i \le f$ and $b_j^+=1+2b_j$,
$1\le j \le n$. Set $A^+=(a_1^+,\ldots,a_f^+)$ and $B^+=(b_1^+,\ldots,b_n^+)$.
Then
$$\aligned
\prod\limits_{i,j} \ (1+a_i+b_j)\ &= 2^{-fn} \prod\limits_{i,j} 
\ [(1+2a_i)+(1+2b_j)] \cr
&= 2^{-fn}\prod\limits_{i,j} (a_i^+ - zb_j^+)|_{z=-1} \cr
&= 2^{-fn} s_{(n)^f}(A^+ - zB^+)|_{z=-1}
\endaligned
$$
where $z$ is a formal $\lambda$-ring element of rank $1$. Here, we have
used the formula (1.1.1) for the resultant.
Moreover, by (1.7.2), 
\smallskip 
$$\aligned
\prod\limits_{i\le j} \ (1 + a_i + a_j) &=
2^{-{{f+1}\choose 2}} \prod\limits_{i\le j} \ [(1+2a_i) + (1+2a_j)] \cr
&=2^{-{f \choose 2}} \ s_{\rho_f}(A^+).
\endaligned
$$
Thus (2.3) can be rewritten as
$$
2^{-N+f} s_{\rho_f}(A^+)\cdot s_{(n)^f}(A^+ - zB^+) |_{z=-1} \ ,
$$
where $N=\rank\,E\vee F$. By the factorization formula (1.6.1)
for elements $A^+$ and $zB^+$,
we can rewrite this last expression as
$$
2^{-N+f}s_{(n+f,n+f-1,\ldots,n+1)}(A^+ - zB^+)|_{z=-1}.
$$
The top-degree component of this polynomial is
$$
2^fs_{(e,e-1,\ldots,n+1)}(A - zB)|_{z=-1}
=2^f\sum_I \ s_{(e,e-1,\ldots,n+1)/I}(A)\cdot s_{\widetilde I}(B),
$$
the sum over partitions $I\subset (e,e-1,\ldots,n+1)$, as desired.
\qed
\enddemo

\head
{\bf 3. The main formulas}
\endhead

The most popular method to compute the fundamental class of 
a subscheme $D\subset X$
tries to find a scheme $G$ mapping properly
to $X$, on which one has a locus $Z$ that maps birationally
onto $D$ and for which one can compute its class $[Z]$.
Usually this is because $[Z]$ is the zero locus of a section
of some bundle whose rank is equal to $\codim_GZ$, so the class $[Z]$
is evaluated to be the top Chern class of the bundle.
For example, this pattern was used in [J-L-P] and many other papers 
(compare [F1] and [F-P]). We will also follow this pattern in the present
section.

\smallskip

We follow the notation from the Introduction and Section 1.
Recall that $n=e-f$ and $q=f-r$.
Let $\pi: G=G^q(F)\to X$ be the Grassmann bundle parametrizing rank $q$
quotients of the bundle $F$. 
On $G$ there exists a tautological sequence
$$
0\to R \to F_{G} \to Q \to 0,
$$
where $\rank R = r$ and $\rank Q = q$. The composite morphism
$$
E_{G}^* @>\varphi_{G}>> F_{G} \toto Q
$$
gives a section of 
$$
H=\Ker (E_{G}\otimes Q \to \wedge^2 Q)
$$
when $\varphi$ is symmetric (and respectively a section of 
$$
H=\Ker (E_{G}\otimes Q \to S^2 Q)
$$ 
when $\varphi$ is skew-symmetric). Let $Z\subset G$ denote the subscheme
of zeros of this section, in both respective cases. Observe that
$\pi$ maps $Z$ onto $D_r(\varphi)$. The next proposition will contain
the main calculation of the paper. 

\smallskip

\proclaim{\bf Proposition 3.1} \ (i) Suppose that, for symmetric $\varphi$,
$\codim_{G}Z=\rank H$ and $\pi$ restricted to $Z$ establishes a birational
isomorphism of $Z$ and $D_r(\varphi)$. Then, the following equality
in $A_*(X)$ holds:
$$
[D_r(\varphi)]= 
\sum_I Q_{\rho_q+I}(F)\cdot s_{\C\widetilde I}(E-F)\cap [X],
\tag 3.2
$$
where the sum is over all partitions $I$ in $(n)^q$ and for the conjugate
partition $\widetilde I\subset (q)^n$ we write $\C\widetilde I 
= (q-\widetilde i_n,\ldots, q-\widetilde i_1)$.

\smallskip

\noindent
(ii) Assume that $r$ is even. Suppose that, for skew-symmetric $\varphi$,  
$\codim_{G}Z=\rank H$ and $\pi$ restricted to $Z$ establishes a birational
isomorphism of $Z$ and $D_r(\varphi)$. Then, 
the following equality in $A_*(X)$ holds:
$$
[D_r(\varphi)]= 
\sum_I P_{\rho_{q-1}+I}(F)\cdot s_{\C\widetilde I}(E-F)\cap [X],
\tag 3.3
$$
where the sum is the same as in (3.2).

\smallskip
\noindent
(iii) Let $n\ge 1$ and assume that $r$ is odd. Suppose that, for  
skew-symmetric $\varphi$, $\codim_GZ=\rank H$ and $\pi$ restricted
to $Z$ establishes a birational isomorphism of $Z$ and $D_r(\varphi)$.
Then, the following equality in $A_*(X)$ holds:
$$
[D_r(\varphi)]= 
\sum_J P_{\rho_q+J}(F)\cdot s_{\C\widetilde J}(E-F)\cap [X],
\tag 3.4
$$
where the sum is over all partitions $J$ in $(n-1)^q$ and for the conjugate
$\widetilde J\subset (q)^{n-1}$ we write $C\widetilde J
=(q-\widetilde j_{n-1},\ldots, q-\widetilde j_1)$.

\endproclaim

\smallskip

\demo{Proof} We give a detailed proof in the symmetric case and make
some necessary comments about the skew-symmetric cases. Let 
$K=\Ker(\alpha: E\toto F)$.
In the Grothendieck group $K(G)$, the following equality holds:
$$
[E\otimes Q]-[\wedge^2Q]=[K_{G}\otimes Q]+[R\otimes Q]+[S^2Q].
$$
Hence we get
$$
[D_r(\varphi)]=\pi_* \Bigl(c_{\top}(K_{G}\otimes Q)\cdot c_{\top}(R\otimes Q) 
\cdot c_{\top}(S^2Q)\cap [G]\Bigr).
\tag 3.5
$$
We have, by (1.7.1),
$$
c_{\top}(K_{G}\otimes Q)=\sum_I s_I(Q)\cdot s_{\C\widetilde I}(K_{G})
=\sum_I s_I(Q)\cdot s_{\C\widetilde I}(E_{G}-F_{G}) ,
$$
where the sum is over partitions $I\subset (n)^q$ and for the conjugate
partition $\widetilde I\subset (q)^n$ \ we write $\C\widetilde I=
(q-\widetilde i_n,\ldots, q-\widetilde i_1)$.

\smallskip

We compute the RHS of (3.5). First, using the above expansion of
$c_{\top}(K_G\otimes Q)$ and (1.7.2), we get
$$
\aligned
&\pi_* \bigl(c_{\top}(K_{G}\otimes Q)\cdot c_{\top}(R\otimes Q) 
\cdot c_{\top}(S^2Q)\cap [G]\bigr) \cr
&=\pi_*\bigl(\sum_I c_{\top}(R\otimes Q)\cdot  Q_{\rho_q}(Q)\cdot s_I(Q) 
\cdot s_{\C\widetilde I}(E_{G}-F_{G})\cap [G]\bigr). 
\endaligned
$$
Secondly, using the factorization formula (1.6.2), we infer that
this last expression equals
$$
\aligned
&\pi_*\bigl(\sum_I c_{\top}(R\otimes Q)\cdot Q_{\rho_q + I}(Q)
\cdot s_{\C\widetilde I}(E_{G}-F_{G})\cap [G]\bigr) \cr
&=\sum_I Q_{\rho_q +I}(F) \cdot s_{\C \widetilde I}(E-F)\cap [X],
\endaligned
$$
where the sum is over $I\subset (n)^q$, and in the last equality
we have used the push-forward formula (1.8.1).

\smallskip

The proof in the skew-symmetric case, when $r$ is even, is analogous; 
in addition, we must
use the push-forward formula (1.8.2).

\smallskip
The proof in the skew-symmetric case, when $n\ge 1$ and $r$ is odd, goes
the same way. We use the push-forward formulas (1.8.3) and (1.8.1).

\smallskip
The proposition has been proved.
\qed
\enddemo

\smallskip
As usual the formulas for the fundamental classes of degeneracy loci
hold under more general assumptions than in the previous proposition.

\proclaim{\bf Theorem 3.6} (i) \ If $X$ is a pure-dimensional Cohen-Macaulay
scheme and $D_r(\varphi)$ is of expected pure codimension $c(r)$ or empty,
then, in the symmetric case, the fundamental class of $D_r(\varphi)$
is evaluated by (3.2).

\smallskip
\noindent
(ii) \ In the skew-symmetric case, under the analogous assumptions,
the fundamental class of $D_r(\varphi)$ is evaluated by (3.3) if $r$ is even,
and by (3.4) if $n\ge 1$ and $r$ is odd.

\endproclaim

\demo{Proof}
We pass to a ``universal case''.(Our method is similar to the technique
explained in [F-P, Appendix A.2].)
For a given morphism $\varphi: E^*\to F$ of one of the two considered types,
we define $\tilde{X}=\Spec \ S^{\sssize\bullet}(E\vee F)^*$ (respectively
$\tilde{X}= \Spec \ S^{\sssize\bullet}(E\wedge F)^*$\,).
Observe that $\varphi$ induces a section $s:X\to \tilde{X}$\,.
On the other hand, there exists the tautological bundle homomorphism 
$\tilde{\varphi}:\tilde E^*\to\tilde F$ 
where $\tilde E=E_{\tilde{X}}$\,,
$\tilde F=F_{\tilde{X}}$\, \ such that
$s^*(\tilde{\varphi})=\varphi$\,.
If $X$ is Cohen-Macaulay, then so is $D_r(\tilde \varphi)$ by virtue
of a result from [DC-L], and because an algebraic fibre bundle
with Cohen-Macaulay base and fibre is Cohen-Macaulay. 
Hence, if $D_r(\varphi)$ is of pure codimension $c(r)$ in $X$, then by 
[F1, Lemma A.7.1] we get 
$$
[D_r(\varphi)]=s^*[D_r(\tilde \varphi)].
\tag 3.7
$$
One checks in a fairly standard way that $\tilde{\varphi}$ 
satisfies the assumptions of Proposition 3.1.
Applying Proposition 3.1 to $\tilde{\varphi}$, we infer that
$[D_r(\tilde \varphi)]$ is evaluated
by (3.2) (resp. by (3.3) or (3.4)) with $\tilde E$ playing the role of $E$ and
$\tilde F$ playing the role of $F$. Consequently, the wanted assertion
follows from (3.7) and the pull-back property of Chern classes.
\qed
\enddemo

\medskip

Note that in the case where $E=F$, Giambelli [G] obtained the expression
$2^qs_{\rho_q}(E)$ for the class of the degeneracy locus considered in
part (i) of the theorem. However, it a classical result of combinatorics
that the following equality of symmetric functions in 
the variables $A=(a_1,a_2,\ldots )$ holds: 
$$
P_{\rho_q}(A)=s_{\rho_q}(A).
$$

\medskip

\example{\bf Example 3.8} \ (i) \ Assume first that $\varphi: E^* \to F$ 
is symmetric.
Let $f=3$ and $r=2$. If $e=4$, then the polynomial representing 
$D_2(\varphi)$ is
$Q_2(F) + Q_1(F)s_1(E-F)$; for $e=5$, the polynomial is
$Q_3(F) + Q_2(F)s_1(E-F) + Q_1(F)s_{1,1}(E-F)$.
\smallskip

For the rest of this example, we assume that $\varphi: E^*\to F$ 
is skew-symmetric.
\smallskip
\noindent
(ii) \ Let $f=3$ and $r=2$. 
If $e=4$ then
the polynomial representing $D_2(\varphi)$ is
$P_1(F) + s_1(E-F)$; for $e=5$, the polynomial is
$P_2(F) + P_1(F)s_1(E-F)+s_{1,1}(E-F)$.
\smallskip
\noindent
(iii) \ If $e=5$, $f=4$,
and $r=2$, then the representing polynomial is
$P_{2,1}(F) + P_2(F)s_1(E-F) + P_1(F)s_2(E-F)$.
\smallskip
\noindent
(iv) \ If $f$ is even and $n=q=1$, then the polynomial representing 
$D_{f-1}(\varphi)$ equals $P_1(F)$.
\smallskip
\noindent
(v) \ If $e=4, f=2$, then $D_1(\varphi)$ is represented by
$P_1(F)s_1(E-F) + P_2(F)$.
\endexample

\smallskip

We leave it to the reader to check that the formulas in this example 
are consistent
with the formulas in (1.11).

\medskip

Perhaps, the easiest way to remember the formula associated with
$(e,f,r)$ in the symmetric case,
is to put $T=(e-r,e-r-1,\ldots, n+1)$, and note
that for $I\subset (e-f)^{f-r}$ the coefficient of $s_{\widetilde I}(E-F)$
is $Q_{T-I}(F)$, $I$ being written {\it increasing} in last index, and
subtraction of the sequences being performed componentwise.
In the skew-symmetric case, with $r$ even, we respectively put $T=(e-r-1,
e-r-2,\ldots, n)$, and note that the coefficient of $s_{\widetilde I}(E-F)$
is $P_{T-I}(F)$ with the same conventions as above. 
A similar interpretation can be given in the skew-symmetric case
when $n\ge 1$ and $r$ is odd.
\smallskip

\example{\bf Example 3.9} \ $e=8, f=4, r=2$; the symmetric case:
$$\aligned
&Q_{6,5}(F) + Q_{6,4}(F)s_1(E-F) + Q_{5,4}(F)s_2(E-F) + Q_{6,3}(F)s_{1,1}(E-F)
\cr 
&+ Q_{5,3}(F)s_{2,1}(E-F)+ Q_{6,2}(F)s_{1,1,1}(E-F) + Q_{5,2}(F)s_{2,1,1}(E-F)
\cr
&+ Q_{4,3}(F)s_{2,2}(E-F) + Q_{6,1}(F)s_{1,1,1,1}(E-F) 
+ Q_{4,2}(F)s_{2,2,1}(E-F)
\cr
&+ Q_{5,1}(F)s_{2,1,1,1}(E-F) + Q_{4,1}(F)s_{2,2,1,1}(E-F) 
+ Q_{3,2}(F)s_{2,2,2}(E-F) \cr 
&+ Q_{3,1}(F)s_{2,2,2,1}(E-F) + Q_{2,1}(F)s_{2,2,2,2}(E-F)
\endaligned
$$
So, in this example, $T=(6,5)$ and, for instance, the coefficient
of $s_{2,2,1,1}(E-F)=s_{\widetilde {(4,2)}}(E-F)$ is $Q_{(6,5)-(2,4)}(F)
=Q_{4,1}(F)$.
\smallskip

A decomposition into a sum of the $s_I(F)\cdot s_J(E)$'s (computed
with the help of the library SFA of ACE (see [V])) \ is $4$ times

$$
\aligned
&s_{6,1}(F)\bigl(s_{1,1,1,1}(E)+s_{2,2}(E)-s_{2,1,1}(E)\bigr)+s_{6,5}(F) \cr
&+s_{6,3}(F)\bigl(-s_2(E)+s_{1,1}(E)\bigr)+s_{2,1}(F)s_{2,2,2,2}(E)
+s_{4,3}(F)s_{2,2}(E)
\cr
&+s_{4,1}(F)\bigl(s_{2,2,1,1}(E)-s_{2,2,2}(E)\bigr)+s_{8,3}(F)+s_{10,1}(F) \cr
&+s_{8,1}(F)\bigl(-s_2(E)+s_{1,1}(E)\bigr).
\endaligned
$$
\endexample

\head
{\bf 4. Some variations}
\endhead 

\smallskip

To compute the fundamental classes of subvarieties, one can also
use appropriate geometric constructions with a {\it nontrivial}
generic fibre.
This method was invented in [P1] in order to give
a short proof of the formulas from [J-L-P] and [H-T], and is summarized
(and somewhat straightened) in the following proposition. 
In this proposition, we may assume
that the Chow groups have rational coefficients. 
  
\smallskip

\proclaim{\bf Proposition 4.1}
Let $D$ be an irreducible (closed) subscheme of a scheme $X$.
Let $\pi:G\to X$ be a proper morphism of schemes and $W$ be a 
(closed) subscheme of 
$G$ such that $\pi(W)=D$.
We have the following two instances:

\noindent
(i) \ Suppose that $G$ is smooth.
Assume that there exists 
$$g\in A_{\dim G + \dim D-\dim W}(G)$$ 
and a point $x$ in the smooth locus of $D$
such that in $A_*(G_{x})$\,, where $G_{x}$ is the fibre of $\pi$ over
$x$, one has: 
$$i^*_x(g)\cdot[W_{x}]=[\point].$$
Here, $W_{x}$ is the fibre of $W$ over $x$ and
$i_x:G_{x} \hoto G$ \ is the
inclusion.
Then the following equality holds in $A_*(X)$\,:
$$
[D]=\pi_*\bigl(g\cdot[W]\bigr)\,.
$$
\smallskip\noindent
(ii) \ Here $G$ is possibly singular. 
Suppose that there exists a family of vector bundles $\{E^{(\alpha)}\}$
on $G$ and $g=P\bigl(\{c.(E^{(\alpha)})\}\bigr)$ a homogeneous 
polynomial 
of degree $\dim W-\dim D$ in the 
Chern classes of $\{E^{(\alpha)}\}$ ($\deg \ c_i(E^{(\alpha)})=i$) 
with rational coefficients,
such that in $A_*(G_{x})$\,,
$$P\bigl(\{c.(i_x^*E^{(\alpha)})\}\bigr)\cap[W_{x}]=[\point],$$
where $x$, $G_x$, $W_{x}$ and $i_{x}$ are as above.
Then the following equality holds in $A_*(X)$\,:
$$
[D]=\pi_*\bigl(g\cap[W]\bigr).
$$
\endproclaim

\demo{Proof}
(i) Using a standard dimension argument, we can replace, in the assertion,
$D$ by its smooth part, i.e., we can assume $D$ is smooth.
Write $G_D=G\times_XD$\,, $W_D=W\times_XD$\,,
$\eta:G_D\to D$ the projection induced by $\pi$, and $k:G_D\to G$
 the inclusion.
Then, the assertion is a consequence of the following identity in $A_*(D)$:
$$
\eta_*\bigl(k^*(g)\cdot[W_D]\bigr)=[D].
$$
To prove this last equation, we first remark that the assumptions imply
$$
\eta_*\bigl(k^*(g)\cdot[W_D]\bigr)=m[D],
$$
where $m\in\Bbb Z$.
Let $x$ be a point in $D$ and consider the fibre square
$$\CD
G_x &\ \ \buildrel j \over{\lhook\joinrel\loto}\ \ & G_D \\
@V p VV  @VV \eta V \\
\{x\} &\buildrel i \over{\lhook\joinrel\loto} &D .
\endCD $$
Using the assumptions on $g$ and [F1, Theorem 6.2], 
we have
$$\aligned
i^*\eta_*\bigl(k^*(g)\cdot[W_D]\bigr)
&= p_*\Bigl(j^*\bigl(k^*(g)\cdot[W_D]\bigr)\Bigr) \\
&= p_*\bigl(i_x^*(g)\cdot[W_x]\bigr)=p_*\bigl([\point]\bigr)=[\point].
\endaligned $$
This implies $m=1$ and assertion (i) is proved.

\bigskip
The proof of (ii) is essentially the same.
\qed
\enddemo
\medskip

Let $F\subset E$ be two vector bundles of ranks $f$ and $e$ on a variety $X$.
We now describe a certain geometric construction associated with 
a bundle morphism $\varphi: F\to E^*$ induced by a section of $E^*\vee F^*$
(resp. $E^*\wedge F^*$). (In this section, we use a slightly different
setup than in the Introduction, Section 1, and Section 3.) 
This construction generalizes in a natural way the construction used
in [P1] and is based on the following characterization of the rank of
the above morphisms.
Assume that $p$ is a natural number such that 
$2p<f$. Let $V\subset U$ be two vector spaces
of dimensions $f$ and $e$. Let $\phi: V \to U^*$ be a linear
map induced by a section of $U^*\vee V^*$ (resp. $U^*\wedge V^*$). 
Then $\rank \phi \le 2p$ iff there
exists a pair of vector spaces $(A,B)$ such that $A$ is a subspace
of $U$ of dimension $e-p$, $B$ is a subspace of $V$ of dimension $f-p$,
$B\subset A$, and the composite map
$$
B \hookrightarrow V @>\ \phi \ >> 
U^* \toto A^*
$$
is zero. Indeed, suppose that such a pair $(A,B)$ exists. Then
the matrix of $\phi$ in some basis has $e-p$ rows whose
initial segments of length $f-p$ consist of zeros. Then every $(2p+1)$-minor
of such a matrix vanishes (use the Laplace expansion of this
minor w.r.t. the first $p+1$ from these $e-p$ rows). The opposite
implication can be showed using fairly standard linear algebra
(by reducing a matrix to its ``standard form" via the elementary
row and column operations).
\medskip

Let
$$
\pi: G=Fl_{f-p,e-p}(F,E)\to X
$$ 
be the flag bundle parametrizing pairs $(A,B)$ where $A$ is a 
rank $e-p$ subbundle of $E$, $B$ is a rank $f-p$ 
subbundle of $F$, and $B\subset A$. Let
$S\subset R$ be the tautological (sub)bundles
of respective ranks $f-p$ and $e-p$ on $G$. 
\smallskip

With $\varphi: F\to E^*$ as above, we associate
a locus $W\subset G$ \ to be the subscheme
of zeros of the composite morphism
$$
S\hookrightarrow F_{G} @>\ \varphi_{G}\ >> 
E^*_{G} \toto R^*.
$$
Let $D=D_{2p}(\varphi)$. By the above discussion, we have $\pi(W)=D$. 

\smallskip

We want now to work with some ``universal" $\varphi$ (like that
in the proof of Theorem 3.6). Moreover, in this situation,
we want now to get information about the generic fibre $W_x=:\Cal F$ of
$\pi\big|_W$ in order to apply Proposition 4.1. To this end,
it suffices to make the following ``universal local study". 
Let $V\subset U$ be vector spaces of dimensions $f$ and $e$, respectively. 
Let $X$ denote the affine space $U^*\vee V^*$ (resp. $U^*\wedge V^*$). 
In this situation, there exists a tautological morphism 
$\varphi: (F=V_X) \to (E=U_X)^*$ and the corresponding subvariety $W$.
Let $\phi: V \to U^*$ be a linear map coming from a 
section of $U^*\vee V^*$ (resp. $U^*\wedge V^*$). 
Then the fibre $W_{\phi}$ over $\phi$ is identified with
$$
W_{\phi} = \{ (L,M) \in Fl \ : \ \ i_L^* \circ \phi \circ j_M = 0 \},
$$
where $Fl=Fl_{f-p,e-p}(V,U)$ and
$i_L: L\hookrightarrow U$, \  $j_M : M\hookrightarrow V$ are the inclusions. 
By calculation in local coordinates, one checks that $W\subset G = Fl$
is a complete intersection of codimension equal to $\rank(R^*\vee S^*)$
(resp. $\rank(R^*\wedge S^*)$). Hence one easily computes that 
$\dim \Cal F = \dim W - \dim D =
p(p-1)/2$ \ (resp. $\dim \Cal F = p(p+1)/2$). Consequently, the dimension of
the generic fibre $\Cal F$ depends only on $p$ and not on $e$ and $f$. 

\smallskip
 
The following very simple fact is helpful to find the class ${g}$ 
satisfying the requirements of Proposition 4.1.

\proclaim{\bf Lemma 4.2} Let $i:Y'\hoto Y$ be a closed embedding of smooth
varieties, let $X\subset Y$ and $X'\subset Y'$ be two subvarieties such that
$i(X')\subset X$ and $\dim\, X'= \dim\, X$\,.
Assume that an element $z\in A^*(X)$ satisfies $[X']\cdot i^*(z)=[\point]$
in $A^*(Y')$.
Then, $[X]\cdot z =[\point]$ in $A^*(Y)$.
\endproclaim
\smallskip

Indeed, we have $i_*[X']=[X]$, and by the projection formula we infer
$$
[\point]=i_*\bigl([X']\cdot i^*(z)\bigr)=i_*[X']\cdot z =[X]\cdot z\,,
$$
as claimed.

\proclaim{\bf Proposition 4.3} \ The class
$g = 2^{-p}s_{\rho_{p-1}}(S^*)$
 \ (resp. $g = s_{\rho_p}(S^*)$) 
satisfies the assumption of Proposition 4.1(ii).
\endproclaim
\smallskip

\demo{Proof}  
We follow the notation from the discussion before the lemma.
Here, the role of ``$x$" from Proposition 4.1 is played by $\phi$
such that $\rank(V @> \phi >> U^* \toto V^*)=2p$.
Let $S'$ denote the tautological rank \ $f-p$ bundle on the Grassmannian
$G_{f-p}(V)$ parametrizing $(f-p)$-dimensional subspaces of $V$.
\smallskip
\noindent
1) \ 
Assume first that $e=f=2p$ \ so $V=U$ and the corresponding bilinear form is
nondegenerate.
Then $[\Cal F]$ is evaluated as the top Chern class of the bundle 
$S^{2}({S'}^{*})$ (resp. ${\wedge^{2}}({S'}^{*})$). 
We get by (1.7),
$$
[\Cal F]=2^ps_{\rho_p}({S'}^*) \qquad 
\hbox{(resp.\ }[\Cal F]=s_{\rho_{p-1}}({S'}^*) \ ).
$$
The assertion now follows by taking the dual Schubert cycles in the
Grassmannian $G^p(V^*)$ (see, e.g., [F1, Chap.14]).
\smallskip
\noindent
2) \
Let now $2p<e=f$ \ (so again $V=U$), and let 
$V'\subset V$ be an inclusion of vector spaces
of dimensions $2p$ and $f$\,, respectively.
Assume that $V$ is endowed with a symmetric (resp. skew-symmetric) 
form $\phi$ of rank $2p$ such that 
the form $\phi \big|_{V'}$ is nondegenerate.
We now use the lemma  with the following data: $Y'=G_p(V')$ and 
$Y=G_{f-p}(V)$\,; $i:G_p(V')\hoto G_{f-p}(V)$ being defined by 
$L\mapsto L\oplus A$\,, where $V=V'\oplus A$\,.
Moreover, $X$ and $X'$ are the generic fibres under consideration and
$z=2^{-p}s_{\rho_{p-1}}({S'}^*)$ (resp. $z=s_{\rho_p}({S'}^*)$ ).
Then part 1) and the lemma yield the desired assertion.
\smallskip
\noindent
3) \ 
Finally, suppose that $f<e$ and let $U=V\oplus B$, where $\dim \ B=n=e-f$.
We now apply the lemma to the following embedding:
$$
i: \bigl(Y'=G_{f-p}(V)\bigr)\hoto \bigl( Y=G \bigr)\,,
$$
where $i(L)=(L,L\!\oplus\!B)$.  
Moreover, $X$ and $X'$ are the generic fibres under consideration and
$z=2^{-p}s_{\rho_{p-1}}(S^*)$ (resp. $z=s_{\rho_p}(S^*)$ ).
Then part 2) and the lemma yield the desired result.
\qed
\enddemo
\bigskip

We will be now interested in formal identities of some polynomials
in Chern classes and their push-forwards. Therefore we assume that
$X$ is smooth and treat the Chern classes as elements of the appropriate
Chow rings.
Using Propositions 4.1(ii) and (4.3), we infer that in the symmetric case 
the degeneracy locus $D$ is represented by 
$$
\pi_* \bigl(c_{\top}(R^*\vee S^*)\cdot 2^{-p} s_{\rho_{p-1}}(S^*) \bigr).
$$
Similarly, in the skew-symmetric case, the degeneracy locus $D$
is represented by
$$
\pi_*\bigl(c_{\top}(R^*\wedge S^*)\cdot s_{\rho_p}(S^* )\bigr).
$$
Of course $D=D_{2p}(\varphi^*)$.
Combining Propositions 2.1 and 2.2 with Theorem 3.6
applied to $\varphi^*: E\to F^*$, we thus
get the following algebraic equalities. In the symmetric case,
writing $T=(e-p,e-p-1,\ldots,n+1)$,
$$
\aligned
&\pi_* \bigl( 2^{-p} \sum_I Q_{\rho_{f-p}+I}(S^*)\cdot 
s_{\C\widetilde I}(R^*-S^*)
\cdot s_{\rho_{p-1}}(S^*)\bigr) \cr
&=\pi_* \bigl(2^{f-2p}\sum_J s_{T/J}(S^*)\cdot
s_{\widetilde J}(R^*-S^*)\cdot s_{\rho_{p-1}}(S^*)\bigr) \cr
&=\sum_L Q_{\rho_{f-2p}+L}(F^*)\cdot s_{\C\widetilde L}(E^*-F^*),
\endaligned
\tag 4.4
$$
where $I$ runs over partitions in $(n)^{f-p}$, $J$ runs over partitions 
in $T$ and $L$ runs over partitions in $(n)^{f-2p}$. 
In the skew-symmetric case, writing $T$ for the partition
$(e-p-1,e-p-2,\ldots,n)$,
$$\aligned
&\pi_* \bigl(\sum_I P_{\rho_{f-p-1}+I}(S^*)\cdot 
s_{\C\widetilde I}(R^*-S^*) \cdot s_{\rho_p}(S^*) \bigr) \cr
&=\pi_* \bigl(\sum_J s_{T/J}(S^*)\cdot s_{\widetilde J}(R^*-S^*)
\cdot s_{\rho_p}(S^*)\bigr) \cr
&=\sum_L P_{\rho_{f-2p-1}+L}(F^*)\cdot s_{\C\widetilde L}(E^*-F^*),
\endaligned
\tag 4.5
$$
where $I$ and $L$ run over the same sets of partitions as above
and $J$ runs over partitions in the present $T$.

\smallskip

These equalities were originally obtained with the help of 
the library SFA of ACE (see [V]) for small values of $f,p$, and $n$. 
We know no algebraic proofs of (4.4) and (4.5), not invoking geometry.
\smallskip
Let $G'=G_{f-p}(F)\to X$ and set 
$$\CD
C=\Coker(S' \hookrightarrow F_{G'} \hookrightarrow E_{G'}),
\endCD
$$
where $S'$ is the tautological subbundle on $G'$.
Then using the presentation
$$
G = G_n(C) \to G_{f-p}(F) \to X
$$
and (1.9), one can rewrite the LHS's of (4.4) and (4.5) as purely
algebraic expressions involving symmetrizing operators. Or, equivalently,
one can embedd $G$ in $G_{f-p}(F) \times_X G_{e-p}(E)$ and
use the product of the symmetrizing operators corresponding to
the factors, at the cost of multiplying the
expressions in (4.4) and (4.5) to be push-forwarded
by $c_{\top}\bigl({S'}^*\otimes (E/R')\bigr)$, where $R'$ is the tautological
subbundle on $G_{e-p}(E)$ (and we omit some pull-back indices).

\bigskip\medskip

Finally, the second author takes this opportunity to make the following

\remark{\bf Remark 4.6}\ (Revisions and corrigenda to [P1], [P4], and [P-R])
\smallskip
\noindent
[P1]: Revisions:  The assumptions on ${g}$ in Proposition 4.1 (in the
present paper) 
straighten an unprecise
expression ``the Poincar\'e dual of'' from $5^0$ in [P1, Proposition 2.1].
The assumption that the ground field is algebraically closed of
characteristic different from $2$ was mistakenly omitted in Section 3.
The reference ``Lemma 9 in [10]'' on p.196 should be replaced by 
``[2, Lemma A.7.1]''.  
\medskip
\noindent
[P4]: Revisions -- insert after ``... following changes:" \ p.$162_5$:
\smallskip
p.154 l. $-10$ \ should read: `` Write $a_i=c_1(L_i^E)$ \ 
$i=1,\ldots,n$. \ ..." \,,
\smallskip
p.154 l.$-5$ should read: `` ... for $\omega$, the fact that $\omega_*$
is induced by $\partial_{(a_2,\ldots,a_n)}^{k-1}$, " ,
\smallskip
p.155 l. $-14$ \ should read: `` ... $= d P_{I,J}(E)\cap 
\alpha$ ..." \,,
\smallskip 
p.155 l. $-13$ \ should read: `` ... Let $I'=(i_2,\ldots,i_k)$." \,, 
\medskip
\noindent
Misprints  -- should be: p.$130_3$ \ ``$(-1)^{p-1}$"  \ //  \ 
p.$130_1$  \ ``$(-1)^p$" \ // \ p.$131^9$ \ ``$\rho_{n-r-1}$" \ // 
 \ p.$148^{11}$ \ ``$-a_{q+1},\ldots $"  \ // \  p.$162_3$ \ 
`` ...$\cdot \xi^{i_1} \cap$ ... " \ // \
p.$164^6$ \ ``$\partial_r \circ \ldots$" \ // \ p.$174^{14}$ \ ` boxes" \ ' .

\medskip
\noindent
[P-R]: Revisions -- should read: p.$19_{18}$ \ ``Whenever, in this paper
we speak about Schubert subschemes in $OG_nV$, we assume that there
exists a completely filtered rank $n$ isotropic subbundle of $V$." \  //
 \ p.$56_{14}$ and p.$60^3$ \ `` ... \ , so we can apply the induction 
assumption to $M_1$ defined below. The partitions ..."  \,. 

\smallskip
\noindent
Misprints -- should be: p.$39_5$ \ `` $\widetilde P_J(X_m)$" \ // \ 
p.$86^{12}$ \ ``Fulton W." \ // \ p.$87_{2}$ \ ``Gieseker-Petri" \,.
 
\endremark

\bigskip\medskip

\noindent
{\bf Acknowledgments.} We thank Bill Fulton for his 
helpful suggestion concerning
Proposition 4.1. We are also grateful to Romek D\c abrowski and 
Witek Kra\'skiewicz for some useful discussions.
The second author thanks the organizers of the
Intersection Theory Conference in Bologna for invitation to attend this 
interesting meeting, as well as for the encouragement to write up
the present paper.

\bigskip\bigskip

\noindent{\eightrm
A.L.: \ C.N.R.S., Institut Gaspard Monge, Universit\'e de Marne-la-Vall\'ee,
5 Bd Descartes, Champs sur Marne, 77454 Marne La Vall\'ee Cedex 2, France

\medskip
\noindent
P.P.: \ Mathematical Institute of Polish Academy of Sciences,
Chopina 12, 87-100 Toru\'n, Poland}

\bye